\documentclass{article}                                                                           

\renewcommand\thesection{\Roman{section}.}                                         
\renewcommand\thesubsection{\thesection\Alph{subsection}.}                   
\renewcommand\thesubsubsection{\thesubsection\arabic{subsubsection}.} 

\usepackage[explicit]{titlesec}                                                                  
\titleformat{\section}{\bfseries}{\thesection}{1em}{\MakeUppercase{#1}}  
\titleformat{\subsubsection}{\itshape}{\thesubsubsection}{1em}{#1}          

\usepackage{setspace}                                                                            
\usepackage{indentfirst}                                                                          
\usepackage[margin=1.3in]{geometry}                                                     

\usepackage{lastpage}                                                                             
\usepackage[figure,table]{totalcount}                                                        

\usepackage{amsmath}                                                                            
\usepackage{multirow}                                                                             
\usepackage{graphicx}                                                                             

\usepackage[]{footmisc}                                                                          

\usepackage{caption}                                                                              
\usepackage[labelformat=simple]{subcaption}                                           
\usepackage{paralist}	
\usepackage{amssymb,amsthm,graphicx}
\usepackage{epsfig}
\usepackage[mathcal]{euscript}
\usepackage{setspace}
\usepackage{color}
\usepackage{array}

\usepackage{float} 
\usepackage{xspace}
\usepackage{mathrsfs}
\usepackage[pdftex]{hyperref}

\newcommand{\Macro}{\ensuremath{\Sigma}}

\newcommand{\ve}[1]{\ensuremath{\mathbf{#1}}}

\begin{document}


\title{Eigenvalue Solvers for Modeling Nuclear Reactors on Leadership Class Machines} 

\author{
\vspace{20mm}
\\R.N.\ Slaybaugh,$^{\text{a},\ast}$ M.\ Ramirez-Zweiger,$^{\text{a}}$ \\Tara Pandya,$^\text{b}$ Steven Hamilton,$^\text{b}$ and T.M. Evans$^\text{b}$\\[4pt] 
\textit{$^a$University of California, Berkeley, Nuclear Engineering Department}\\[-10pt]       
\textit{4173 Etcheverry Hall, Berkeley, CA 94720, USA} \\[-5pt]
\textit{$^b$Oak Ridge National Laboratory, Radiation Transport and Criticality Group} \\ [-10pt]
\textit{P.O. Box 2008, Oak Ridge, TN 37831-6170, USA} \\ [-2pt]
{$^\ast$slaybaugh@berkeley.edu}}       

\date{                               
\vspace{40mm}
Number of pages: \pageref{LastPage} \\  
Number of tables: \totaltables \\
Number of figures: \totalfigures \\}                                                                                           

\maketitle

\pagebreak

\begin{abstract}
{Three complementary methods have been implemented in the code Denovo  that accelerate neutral particle transport calculations with methods that use leadership-class computers fully and effectively: a multigroup block (MG) Krylov solver, a Rayleigh Quotient Iteration (RQI) eigenvalue solver, and a multigrid in energy (MGE) preconditioner. The MG Krylov solver converges more quickly than Gauss Seidel and enables energy decomposition such that Denovo can scale to hundreds of thousands of cores. RQI should converge in fewer iterations than power iteration (PI) for large and challenging problems. RQI creates shifted systems that would not be tractable without the MG Krylov solver. It also creates ill-conditioned matrices. The MGE preconditioner reduces iteration count significantly when used with RQI and takes advantage of the new energy decomposition such that it can scale efficiently. Each individual method has been described before, but this is the first time they have been demonstrated to work together effectively. 

The combination of solvers enables the RQI eigenvalue solver to work better than the other available solvers for large reactors problems on leadership class machines. Using these methods together, RQI converged in fewer iterations and in less time than PI for a full pressurized water reactor core. These solvers also performed better than an Arnoldi eigenvalue solver for a reactor benchmark problem when energy decomposition is needed. The MG Krylov, MGE preconditioner, and RQI solver combination also scales well in energy. This solver set is a strong choice for very large and challenging problems. 

}
\end{abstract}

\pagebreak


\section{Introduction}
\label{sec:intro}
The steady-state Boltzmann equation for neutron transport covers six dimensions of phase space. Typical deterministic transport problems today are three-dimensional, have up to thousands $\times$ thousands $\times$ thousands of mesh points, use up to $\sim$150 energy groups, include accurate expansions of scattering terms, and are solved over many angular directions. The next generation of challenging problems are even more highly refined. High-fidelity, coupled, multiphysics calculations are the new ``grand challenge'' problems for reactor analysis, requiring that the finely-resolved neutron flux be calculated quickly and accurately.

Very large computers, such as Titan \cite{Titan2013}, are available to perform such high-fidelity calculations. Historical solution methods are not able to take full advantage of new computer architectures, or they have convergence properties that limit their usefulness for difficult problems. The goal of this research is to accelerate transport calculations with methods that use new computers fully and effectively, facilitating the design of better nuclear systems. 

Three complimentary methods have been implemented in the code Denovo \cite{Evans2010} that accomplish this goal: a multigroup block (MG) Krylov solver, a Rayleigh Quotient Iteration (RQI) eigenvalue solver, and a multigrid in energy (MGE) preconditioner. Each individual method has been generally described before (see \cite{Slaybaugh2012}, \cite{Slaybaugh2013}), but this is the first time they have been demonstrated to work together in a complementary way. 

The driving concept of this research is to use RQI for solving the $k$-eigenvalue problem in 3-D neutron transport. This objective was not tractable with the tools originally available in Denovo. The MG Krylov solver and MGE preconditioner were developed to facilitate RQI, though both of these tools are also useful on their own. 

The MG Krylov solver was designed to improve convergence when compared to Gauss Seidel (GS) and to dramatically increase the number cores Denovo can use. Instead of sequentially solving each group with some inner iteration method and then using GS for outer iterations to converge the upscattering, the MG Krylov solver treats a block of groups (either all groups or just upscattering groups) at once such that the inner-outer iteration structure is removed. This results in faster convergence for most problem types. In addition, the block Krylov solver allows energy groups to be solved simultaneously because the multigroup-sized matrix vector multiply can be divided up in energy and parallelized. This extends the number of cores that can be used efficiently by Denovo from tens of thousands to hundreds of thousands \cite{Davidson2013}.

A MGE preconditioner was added to Denovo to reduce iteration count for all problem types, and to address convergence issues associated with RQI. The  preconditioner conducts a multigrid method in the energy dimension. A set of energy grids with increasingly coarse energy group structures are created. This is implemented in a way that easily and efficiently takes advantage of the new energy decomposition. The multigrid algorithm is applied within each energy set such that the energy groups are only restricted and prolonged between groups on that set. Sets do not communicate with one another in the preconditioner, so the scaling in energy is very good. 

Theory indicates that RQI should converge in fewer iterations than traditional eigenvalue solvers like Power Iteration (PI), particularly for problems that are challenging for those solvers. However, the use of RQI would not be practical without the MG Krylov solver and the MGE preconditioner. The implementation of RQI results in a set of equations that is mathematically equivalent to having upscattering in every group, so the scattering matrix becomes energy-block dense. Handling energy-block dense systems when there are more than a few energy groups is not tractable with GS as the multigroup solver. It is only the MG Krylov solver that makes RQI reasonable to use when there are many energy groups. In addition, the MGE preconditioner is needed to mitigate the slow convergence associated with Krylov methods when trying to solve the ill-conditioned systems created by RQI. 

The remainder of this paper presents information about why these methods are complementary as well as results demonstrating that they are. Section \ref{sec:background} discusses each of the new methods in the context of commonly-used methods. Section \ref{sec:pastwork} gives an overview of relevant past work. New results from using the three new methods together are shown in Section \ref{sec:results}, and concluding remarks are made in Section \ref{sec:conclusions}.

\section{Background}
\label{sec:background}
The steady state Boltzmann transport equation, discretized in energy (multigroup), space, and angle ($S_N$) can be written in operator form as
\begin{alignat}{2}
  \ve{L}\psi &= \ve{MS}\phi + q \:, \qquad &\text{(fixed source)} \label{eq:fxdsource} \\
  \ve{L}\psi &= \ve{MS}\phi + \frac{1}{k}\ve{M}\chi f^{T}\phi \:. \qquad &\text{(eigenvalue)} 
  \label{eq:eigenvalue}
\end{alignat}
Here, $\ve{L}$ is the first-order linear differential transport operator; $\ve{M}$ is the moment-to-discrete operator that projects the angular flux moments, $\phi$, onto discrete angles; $\ve{S}$ is the scattering matrix; $q$ is a source term; $f$ is the fission operator, $\nu \Macro_{f}$; $\chi$ is the energy distribution with which neutrons are born out of fission; and $k$ is the asymptotic ratio of the number of neutrons in one generation to the number in the next. The angular flux moments are related to the angular flux, $\psi$, through the discrete-to-moment operator: $\phi = \mathbf{D} \psi$. Using this relationship, Equations \eqref{eq:fxdsource} and \eqref{eq:eigenvalue} can be rearranged such that they are a function of only $\phi$. The formulation is aided by defining $\ve{T} = \ve{DL}^{-1}$ and $\ve{F} = \chi f^{T}$ \cite{Evans2011}:
\begin{alignat}{2}
  (\ve{I} - \ve{TMS})\phi &= q \:, \qquad &\text{(fixed source)} \label{eq:OperatorFxdForm} \\
  (\ve{I} - \ve{TMS})\phi &= \frac{1}{k} \ve{TMF} \phi \:. &\text{(eigenvalue)} \label{eq:OperatorEvalForm}
\end{alignat}

Once the matrices are multiplied together, a series of single ``within-group'' equations that are each only a function of space and angle result. If the groups are coupled together by neutrons scattering from a low energy group to a higher energy group (upscattering), then iterative ``multigroup'' solves over the coupled portion of the energy range may be required. If the eigenvalue is desired, an additional ``eigenvalue'' solve is needed, where $k$ is the dominant eigenvalue and $\phi$ is the corresponding eigenvector \cite{Evans2009}.

\subsection{Block Krylov Solver}
\label{sec:blockkrylov}
Traditionally, the multigroup solve has been done with GS, which is iterative in energy. A space-angle solve using a within-group solver, such as source iteration or a Krylov method, is performed for each energy group in series. The groups are solved from $g=0$, the highest energy, to $g=G$, the lowest. For a group $g$ and an energy iteration index $j$ this is \cite{Evans2010}
\begin{equation}
  \bigl( \ve{I} - \ve{TMS}_{gg} \bigr) \phi^{j+1}_{g} = \ve{TM} \bigl( \sum_{g'=0}^{g-1}\ve{S}_{gg'}\phi^{j+1}_{g'} + \sum_{g'=g+1}^{G} \ve{S}_{gg'}\phi^{j}_{g'}  + q_{g} \bigr)  \:.
 \label{eq:up-GS}
\end{equation}

The first term on the right includes downscattering contributions from higher energies, and the second term represents upscattering contributions from lower energy groups that have not yet been converged for this energy iteration. Groups that only contain downscattering are simply solved once since the second term on the right is zero. Groups with upscattering, however, must be iterated until they converge. Convergence of GS is governed by the spectral radius of the system, so the method can be very slow when upscattering has a large influence on the solution \cite{Adams2002}. GS is fundamentally serial in energy because of how the group-to-group coupling is treated. 

The MG Krylov solver removes the traditional ``within-group'' / ``multigroup'' iteration structure by combining the space-angle and energy iterations to make one space-angle-energy iteration level. This allows the energy groups to be decomposed such that they can be solved in parallel. The space-angle-energy iterations are much like the within-group space-angle iterations, except that the iteration is over a block of groups instead of just one group.  In Denovo, energy may be decomposed over all of the groups (termed full partitioning) or only over the upscatter (thermal) groups (termed upscatter partitioning). With upscatter partitioning, the downscatter groups are replicated and solved using GS while the upscatter groups are solved with a Krylov solver. With full partitioning, only the Krylov solver is used.

The MG Krylov method applied to the upscattering block is shown here, where $\ve{S}_{\text{up\_block}}$ contains the upscattering groups and $\ve{S}_{\text{up\_source}}$ has the downscattering-only groups and $n$ is the space-angle-energy iteration index:
\begin{equation}
  \underbrace{(\ve{I} - \ve{TMS}_{\text{up\_block}})}_{\tilde{\ve{A}}}\phi_{\text{up\_block}}^{n+1} = \ve{TM}(\ve{S}_{\text{up\_source}}\phi_{\text{up\_source}}^{n+1} + q) \:.
  \label{eq:MGkrylov}
\end{equation}
Trilinos \cite{1089021} provides Denovo's Krylov solver, with a choice of either GMRES($m$) or BiCGSTAB \cite{Evans2010}. The Krylov solver is given an operator that implements the action of $\ve{\tilde{A}}$, or the matrix-vector multiply and sweep. In the MG Krylov solver, $\ve{\tilde{A}}$ is applied to an iteration vector, $v$, containing the entire upscattering block instead of just one group.
%

To implement the energy parallelization, the problem is divided into energy sets, with groups distributed evenly among sets. After each set performs its part of the matrix-vector multiply, a global reduce-plus-scatter is the only required inter-set communication. Since each set uses the entire spatial mesh with the same spatial decomposition, the established performance of spatial scaling does not change. The space-angle decomposition in Denovo comes from the KBA wavefront algorithm \cite{Baker1998}. 

The added energy decomposition offers the ability to further decompose a problem, even if the performance limit of spatial decomposition has been reached. The total number of cores is equal to the number of computational domains, that is, the product of the number of energy sets and the number of spatial blocks. For example, a problem decomposed into 20,000 spatial blocks and 10 energy sets would use 200,000 cores. See Ref.\ \cite{Evans2011} for more details.

The MG Krylov solver has been shown to successfully scale to hundreds of thousands of cores. For example, a fixed-source test scaled from 69,102 cores to 190,080 cores with 98\% efficiency \cite{Slaybaugh2011}. Upscatter partitioning has been shown to initially perform better than full partitioning, but does not scale well over many energy sets.  An added benefit of this solver is that Krylov methods generally converge more quickly than GS for problems with upscattering \cite{Davidson2013}.

\subsection{Eigenvalue Solvers}
\label{sec:eigenvalue}
Denovo has three eigenvalue solver choices, PI, RQI, and Arnoldi. 

\subsubsection{Power Iteration}
A common way to solve $k$-eigenvalue problems is with PI. This method is attractive because it only requires matrix-vector products and two vectors of storage space. 
\begin{align}
  \ve{A}\phi &= k\phi \:, \label{eq:EnergyDepEval} \\
  \qquad \text{where}  \qquad \ve{A} &= (\ve{I} - \ve{TMS})^{-1} \ve{TMF} \:, \nonumber \\
  \phi^{i+1} = \frac{1}{k^i}\ve{A}\phi^{i} &\:; \qquad 
  k^{i+1} = k^i \frac{\Vert f^T \phi^{i+1} \Vert}{\Vert f^T \phi^i\Vert} \:.
  \label{eq:PowerIteration}
\end{align} 
PI uses the form of the problem seen in Eq. \eqref{eq:EnergyDepEval} and then iterates as shown in Eq.\ \eqref{eq:PowerIteration}, where $i$ is the iteration index. This converts the generalized form of the eigenvalue problem seen in Eq.\ \eqref{eq:OperatorEvalForm} to the standard form. In the generalized form, the eigenvector-value pair is $(\phi, \frac{1}{k})$, and in the standard form it is $(\phi, k)$.
In legacy applications, the eigenvector is often the fission source rather than the flux moments \cite{Lewis1993,Evans2011}.
Inside of PI, the application of $\ve{A}$ to $\phi$ requires the solution of a multigroup linear system that looks like a fixed source problem,
\begin{equation}
  (\ve{I} - \ve{TMS})\phi^{i+1} = \frac{1}{k^i}\ve{TMF}\phi^{i} \:. \label{eq:EvalDepFxdSource}
\end{equation}

PI's convergence can be very slow for problems of interest.  The error from PI is reduced in each iteration by a factor of $\ve{A}$'s dominance ratio, $\frac{\lambda_{2}}{\lambda_{1}}$. For large, loosely coupled systems, $\lambda_2 \approx \lambda_1$ and PI will therefore converge slowly.

\subsubsection{Rayleigh Quotient Iteration}
Shifted inverse iteration (SII) typically converges more quickly than PI. SII capitalizes on the fact that for some shift $\mu$, $(\ve{A} - \mu \ve{I})$ will have the same eigenvectors as $\ve{A}$ and eigenvalues of $\ve{A}$ that are near the shift will be transformed to extremal eigenvalues of $(\ve{A} - \mu \ve{I})$ that are well separated from the others. The shifted and inverted matrix is used in a PI-type scheme. Given a good shift, $\mu \approx \lambda_1$, SII usually converges more quickly than PI, especially for loosely coupled systems \cite{Allen2002}.


RQI is a shifted inverse iteration method
that uses a dynamically updated shift: the Rayleigh quotient (RQ). For a generalized eigenvalue problem $\ve{A}x = \lambda \ve{B}x$, the RQ is given by
\begin{equation}
  \rho = \frac{x^{T} \ve{A} x}{x^{T} \ve{B} x} \:.
  \label{eq:RQ}
\end{equation}
If $x$ is an eigenvector corresponding to $\lambda$, then $\rho \equiv \lambda$.
The idea of RQI is to use the RQ as the shift in SII using the current eigenvector estimate:
\begin{equation}
  \mu_i = \frac{x_i^{T} \ve{A} x_i}{x_i^{T} \ve{B} x_i} \:.
  \label{eq:Shift}
\end{equation} 
For symmetric systems, the use of the RQ as the shift is optimal and leads to cubic convergence to an eigenvalue/eigenvector pair under appropriate assumptions \cite{Parlett1974}.  The theory for the nonsymmetric case is less rigorous, but nevertheless displays locally quadratic convergence under suitable assumptions \cite{SaadEig}.

RQI has been implemented in Denovo, as detailed in Ref.\ \cite{Slaybaugh2012}, by subtracting $\rho \ve{TMF}$ from both sides of Eq.\ \eqref{eq:OperatorEvalForm}. This gives the following shifted system, where $\gamma \equiv \frac{1}{k}$:
\begin{align}
  (\ve{I} - \ve{TM}\ve{\tilde{S}})\phi &=( \gamma - \rho) \ve{TMF} \phi  \:, 
  \label{eq:OperatorShiftedEval} \\
  \text{ where } \ve{\tilde{S}} &\equiv \ve{S} + \rho\ve{F}  \nonumber \:.
\end{align}
The scalar $(\gamma - \rho)$ on the right hand side may introduce numerical scaling issues when $\gamma \approx \rho$ and can be omitted without altering the convergence behavior of the method.  While $\ve{S}$ is predominantly lower triangular, the addition of $\ve{F}$ makes $\ve{\tilde{S}}$ essentially dense in energy.
Traditional solution methods for the fixed source part of the equation do not handle dense scattering matrices well. In fact, GS may fail to converge when $\ve{F}$ is included.  This has hampered the implementation of SII in multi-group, 3-D codes.

The RQI method added to Denovo uses the MG Krylov solver, which is designed to handle dense scattering matrices, effectively overcoming the burden of using of a dense $\ve{\tilde{S}}$. In addition, RQI can be decomposed in energy and take advantage of the scaling properties of the multigroup solver.


On initial consideration, it would seem shifted inverse methods might not work well when the shift is very good because the matrix becomes so ill-conditioned. Peters and Wilkinson \cite{Peters1979}, however, proved that ill-conditioning is not a fundamental problem for inverse iteration methods. Trefethen and Bau \cite{Trefethen1997} assert that this is the case as long as the fixed source portion is solved with a backwards stable algorithm. 
Paige et al.\ \cite{Paige2006} demonstrated that GMRES is backwards stable when finding $x$ in $\ve{A}x = b$ for a ``sufficiently nonsingular $\ve{A}$'', and define associated criteria. Many researchers have found that Krylov methods must be preconditioned to be able to get good results in practice \cite{Benzi2002}, \cite{Saad1986}, \cite{Trefethen1997} , \cite{Paige2006}. The lack of an effective preconditioner was a limiting factor in the approach described in Ref.~\cite{Slaybaugh2012}.

\subsubsection{Arnoldi}
Denovo also has an Arnoldi Krylov subspace solver available for solving eigenvalue problems, which can take an energy-dependent or energy-independent form \cite{Davidson2013}. Recall the energy dependent form of the eigenvalue equation, as seen in Eqn.~\eqref{eq:EnergyDepEval}. 
%
At each iteration, we apply the operator in Eqn.~\eqref{eq:EnergyDepEval} by first
computing a right hand side, with iteration index $h$, as
\begin{equation}
z^{(h)} = \mathbf{TMF}\nu^{(h)}\:,
\end{equation}
and then solving the linear system
\begin{equation}
(\mathbf{I}- \mathbf{TMS})y^{(h)} = z^{(h)}\:,
\end{equation}
where $z$ is a multi-group-sized vector like in RQI (recall that this vector
covers only one group in PI). Using the multi-group-sized vector allows us to
take advantage of all of the machinery of the block Krylov solver.
 
The energy independent case takes the form
\begin{equation}
\begin{split}
\mathbf{A}\Gamma &= k\Gamma, \\
\mathbf{A} &= \mathbf{f}^T(\mathbf{I - TMS})^{-1}\mathbf{TM\chi}.
\end{split}
\end{equation}
Like in the energy-dependent case, a matrix-vector multiply and sweep and a fixed-source solve are performed, but now there is an additional matrix-vector multiply afterwards:
\begin{align}
z^{(h)} &= \mathbf{TM\chi}\nu^{(h)}\:,\\
(\mathbf{I} &- \mathbf{TMS})x^{(h)} = z^{(h)}\:,\\
y^{(h)} &= \mathbf{f}^Tx^{(h)}\:.
\end{align}

Note that the vectors in the energy-dependent subspace span all of the energy groups whereas the energy-independent approach requires each vector to span just one group--using less memory. However, when the flux moments as well as the eigenvalue are desired, the energy-independent approach requires a final fixed source solve, potentially losing its memory advantage. The fixed-source solve is given by
\begin{equation}
(\mathbf{I} - \mathbf{TMS})\phi = \frac{1}{k}\mathbf{TM\chi\Gamma},
\end{equation}
where $k$ is the eigenvalue and $\Gamma$ is the eigenvector. Note that the nonzero eigenvalues of the energy-dependent and energy-independent formulations are identical; therefore, the convergence behavior of the two methods is identical.  We use the energy-dependent implementation in all of our calculations. As with RQI, the Arnoldi solver can be parallelized over energy by using energy sets enabled by the MG Krylov solver.

\subsection{Multigrid in Energy Preconditioner}
\label{sec:precond}
Preconditioning is important for increasing the robustness of Krylov methods and decreasing Krylov iteration count. This is particularly true for the MG Krylov solver. This solver can create large Krylov subspaces because it forms the subspaces with multiple-group-sized vectors. 
Each additional Krylov iteration increases the subspace size by one multi-group-sized vector and the cost of an iteration application correspondingly increases. 
Therefore, any reduction in iteration count will have a significant benefit in terms of memory and cost per iteration. 
 
Right preconditioning leaves the right hand side of the equation unaffected and does not change the norm of the residual, which is used for convergence testing in most iterative methods. A right preconditioner that does multigrid in the energy dimension and is designed to work with the MG Krylov solver was implemented in Denovo \cite{Slaybaugh2013}. To understand why MGE makes sense for neutron transport, some highlights about these methods are discussed here (derived from Ref.\ \cite{Briggs2000}). 

The error in $x_i$, the $i$th guess for $\ve{A}x_i=b_i$, can be written as a combination of Fourier modes. Each Fourier mode has a frequency, and the frequencies can range from low-frequency (smooth) to high-frequency (oscillatory). Iterative methods, also referred to as smoothers or relaxers, remove high-frequency error components quickly, but take many iterations to remove the low-frequency ones. 

The idea of multigrid methods is to take advantage of the smoothing effects of iterative methods by making smooth errors look oscillatory and thus easier to remove. Errors that are low-frequency on a fine grid can be mapped onto a coarser grid where they are high-frequency. A relaxer is applied on the coarser grid to remove the now oscillatory error components. The remaining error is mapped to a still coarser grid and smoothed again. The problem is restricted to coarser and coarser grids until the coarsest grid is reached.

Next, the coarsest result is prolonged back to the next-finer grid and used to correct the solution there. A few relaxations are done on this finer grid. The errors are prolonged back up the grid structure, continuously correcting on finer grids, until the finest grid is reached. This entire process is called a V-cycle. 
 
Thus, multigrid methods remove the low-frequency error modes that require many Krylov iterations. The preconditioner was designed to take advantage of the energy decomposition used by the MG Krylov method. Each energy set does work only on its own grids and does not need to communicate with other energy sets. This is a communication savings compared to using grids in space or angle. An additional benefit is the simplicity of energy grids. Energy is one-dimensional, which allows for simpler coarsening and refinement than spatial or angular grids.

To implement MGE as a right preconditioner in Denovo, the problem is broken into two steps shown below. Here, $\ve{G}^{-1}$ represents the application of the preconditioner and $y$ is defined as $\ve{G}\phi$; recall that $\ve{A} = \ve{I} - \ve{TMS}$.  
\begin{enumerate}
  \item With a Krylov method solve 
    \begin{equation}
      \ve{AG}^{-1}y = b \:; \label{eq:PrecondKrylov} 
    \end{equation}
  \item after finding $y$, calculate 
    \begin{equation}
      \phi = \ve{G}^{-1}y \:. \label{eq:PrecondPhi}
    \end{equation}
\end{enumerate}

Note: it is possible to implement the preconditioner using the shifted operator (used in RQI) rather than the unshifted operator. In that case $\mathbf{S}$ becomes $\tilde{\ve{S}} = \ve{S} + \rho\ve{F}$ and the right hand side operator becomes $(\frac{1}{k} - \rho)\ve{TMF}$. In practice and as expected, the unshifted operator is a much better choice because it is more stable. 

In this preconditioner the grids are in energy, where the energy group structure is coarsened so that each lower grid has fewer groups. The finest grid is the input energy structure, and the coarsest grid has one or a few groups. Each level has half as many groups as the previous level, rounded up if applicable. If there are $G+1$ groups on the fine grid there will be either $\frac{G+1}{2}$ or $\frac{G+2}{2}$ groups on the coarse grid. This is conceptually straightforward because the energy groups can be combined (restricted) and separated (prolonged) linearly. For details on the implementation of the preconditioner consult Ref. \cite{Slaybaugh2013}
%

The user chooses the number of V-cycles done for each preconditioner application. One V-cycle proceeds from the finest grid to the coarsest grid and back to the finest. Each additional V-cycle should remove more error, but has a computational cost. The depth of the V-cycle can also be specified by the user. The default behavior is determined by the number of groups, such that the grids will be coarsened until there is only one energy group. The number of grids needed is $\text{floor}\bigl( \log_{2}(G-1) \bigr) + 2$ \cite{BinaryTree2012}.

When using multiple energy sets, each energy set does its own ``mini'' V-cycle. 
Each set restricts, prolongs, and relaxes on only its own groups. 
This means that cross-set communication is not needed, which is advantageous for scaling in energy.
In the case of an unequal number of groups per set, all sets use the shallowest grid depth, which is determined by the set with the fewest groups, to enforce energy load balancing between sets. Thus, each set restricts to one or two groups, giving approximately $num\_sets$ total groups across all sets at the coarsest level. The number of grids needed is determined by the set with the minimum number of groups, since it will be the first to reach a grid with one group:
\begin{align}
  num\_g_{min} &= \text{floor}\bigl(\frac{num\_groups}{num\_sets}\bigr) \:, \\
  num\_grids &= \text{floor}\bigl( \log_{2}(num\_g_{min}) \bigr) + 2 \:.
  \label{eq:multisetGrids}
\end{align}

Some number of relaxations are performed on each level while traversing down and up the grids in a V-cycle. Performing more relaxations per grid should remove more error, but has a computational cost. The implemented relaxation method is weighted Richardson Iteration. When applied to the transport equation, this is
\begin{equation}
  \phi^{m} = \bigr(\ve{I} + \omega(\ve{TMS} - \ve{I})\bigl)\phi^{m-1} + \omega b^{m-1} \:,
  \label{eq:relax}
 \end{equation}
where $\omega$ is a constant selected by the user that defaults to 1 and $m$ is the iteration index. 

An important principle is that the preconditioner is only attempting to approximately invert $\ve{A}$. It is therefore reasonable to use a less accurate angular discretization in the preconditioner than the rest of the code. For example, the whole problem may be solved at $S_{10}$, but the preconditioner could use $S_{2}$. The user can specify an angular quadrature set to use in the preconditioner that is different from the angular quadrature set used in the rest of the problem; the default is to use the same quadrature in both. At this time, this option has only been implemented for vacuum boundary conditions. 

\section{Past Work}
\label{sec:pastwork}
This section is a summary of results using RQI without preconditioning \cite{Slaybaugh2012} and the MGE preconditioner with fixed source problems or PI \cite{Slaybaugh2013}. The purpose of this section is to highlight the capabilities that have been demonstrated and point out the short-comings that can be overcome by using the MG Krylov solver, RQI eigenvalue solver, and MGE preconditioner together. In the following discussion, reducing the Krylov iteration count is the primary measure of success--with reduced time as the second consideration--as the software was not implemented optimally at the time of these studies. 

\subsection{Unpreconditioned RQI}
The goal of the RQI studies that have been published was to find out it if RQI useful without preconditioning. We solved two small eigenvalue test problems, one with vacuum boundaries and one with reflecting, that had small dominance ratios. We found that RQI got the correct answer and converged in fewer iterations than PI. However, an intermediately-sized problem did not work so well. Each multigroup solve, except the first, used the maximum number of Krylov multigroup iterations per eigenvalue iteration. This means the eigenvector did not converge after the first iteration. The value of $k$ oscillated between 0.3966 and 0.3967 (the correct value was 0.4) until the calculation was manually terminated. In two more realistic calculations, the 2D and 3D C5G7 MOX Benchmark problems (\cite{OECD-NEA2003}, \cite{OECD-NEA2005}), RQI did not converge the eigenvector nor find an eigenvalue close to the correct one. 

These more challenging problems showed that, as expected, the Krylov solver often cannot converge the eigenvector with the ill-conditioned systems created by RQI. When the Krylov iterations do not converge, the flux estimate is not good. Without an adequate approximation to the eigenvector, the RQ is no longer a valid approximation to the eigenvalue, and therefore the eigenvalue problem does not converge. 

The small problems showed that RQI can require fewer Krylov iterations than PI, and has the potential to be beneficial if the multigroup iterations are converged. If the MG Krylov solver is preconditioned so that the eigenvector converges, RQI may be able to find the correct eigenvalue more efficiently than PI for cases of interest. This leads to the question: what preconditioner should we choose?

Iterative methods reduce oscillatory error modes effectively, but not smooth error modes. The smooth error can prevent iterative methods from converging. This behavior is characterized by rapid error reduction in the first several iterations followed by very little error reduction. Such a trend was observed in tests where RQI failed. Multigrid methods selectively remove smoother error components and are therefore strong candidates for improving convergence in this type of problem. Thus, a multigrid preconditioner should work very well with RQI. 

\subsection{MGE Preconditioner}
Much of the previously published work for the MGE preconditioner focused on choosing all of the options that control the preconditioner. The preconditioning parameters are the Richardson iteration weight, $w$, the number of V-cycles per preconditioner application, and the number of relaxations per level. The syntax used throughout this work will be that $w\#$ is the weight, $r\#$ is the number of relaxations per level, and $v\#$ is the number of V-cycles, e.g.\ $w1r1v1$ is one relaxation per level, one V-cycle, and a weight of one. Using more preconditioning means using larger values of $w$ and/or $r$ and/or $v$. 

Other issues investigated were using a different quadrature set inside the preconditioner than in the rest of the problem, strong scaling, and changing the depth of the V-cycle. 
Tests showed that using a reduced angle set inside the preconditioner is very valuable. For example, using $S_2$ inside MGE for an iron-graphite test solved with $S_8$ reduced the solve time by 73\% compared to using $S_8$ in the preconditioner.

Another important area of investigation was how the preconditioner fared when using multiple energy sets because MGE was built to take advantage of multiple energy sets. The problems solved scaled about 40\% more efficiently with MGE preconditioning than without it. The main reason for the preconditioner's good scaling is that as the number of sets increases, each application of the preconditioner becomes less costly. 

With multisets, the total preconditioning cost goes down because the V-cycle becomes shallower. When 27 groups are used with one set, six grids are created. When 27 groups are used with 10 sets, two grids are created for each set. With fewer grids, each application of the preconditioner performs fewer total relaxations, and is therefore less time intensive. The number of GMRES iterations did not change with the number of sets in the preconditioned cases. This means convergence improvement from the preconditioner did not come from the depth of the V-cycle, at least not for the problems tested.

The results of the multiple energy set study prompted an investigation of controlling V-cycle depth explicitly. The results from several tests confirm the multiple energy set study findings: using only a few grids is better than using many. The optimal number of grids will be problem dependent, but a default grid depth of two was recommended.

All of these tests inform the best way to use the MGE preconditioner, but do not determine whether and when it is useful. To begin investigating this question, two eigenvalue problems were solved with preconditioned and unpreconditioned PI. The first problem was the 3D C5G7 MOX Benchmark problem. In this case the use of MGE significantly reduced Krylov iteration count, but increased the overall runtime. A full PWR problem (described in subsection \ref{subsec:PWR}) exhibited similar behavior: fewer iterations in more time. 

These results suggest that the MGE preconditioner was a failed experiment after all. However, the two eigenvalue problems solved were not particularly challenging for PI, meaning they may not have needed preconditioning to begin with. More importantly, the mathematical properties of the MGE preconditioner suggest it would benefit the RQI solver. Thus, the work published so far has not settled the question of whether MGE is a useful preconditioner for at least some problems. 

\section{Results}
\label{sec:results}
The collection of observations in the past work section led to the questions: will preconditioning with MGE facilitate the use of RQI, and will the combination of RQI, MGE, and the block Krylov solver be advantageous for at least some problems of interest? This section covers a series of problems designed to answer these questions.

We will start by using a relatively small problem to perform a parameter study to investigate the best choice for $r\#v\#w\#$ when using MGE with RQI.
We then solve two much larger problems on a few different machines to answer the posed questions.

Unless otherwise noted, all test problems used a step characteristic (SC) spatial solver (1 unknown/cell and positive in the presence of physical scattering\footnote{Step-characteristics are only guaranteed to be positive in the presence of positive scattering sources. Negativities can still occur when the source is negative due to truncation resulting from $P_n$ anisotropic scattering. Negativities can also result from incomplete convergence in GMRES.}), level-symmetric angular quadrature, the initial guess shift for RQI is 1.0, and the MG grid depth was determined using the default approach. 
The Krylov solver was GMRES, which is set to limit the number of multigroup iterations to 1,000 if the problem does not converge earlier. The convergence tolerances are noted for each problem. The tolerance for the multigroup solve is the convergence tolerance used by GMRES in Trilinos \cite{1089021}. The eigenvalue tolerance is used by the eigenvalue solver to determine if the eigenvalue has converged. In Denovo, PI also checks the L2-norm and the infinity-norm of the difference in the fission source between iterations. The default L2-norm tolerance is 1.0 and infinity-norm tolerance is 0.01.

\subsection{MGE Parameter Choice}
\label{subsec:mge-params}

Since the previous work investigating the best MGE parameter choice, many of the internal workings of Denovo have been updated. 
Further, most problems studied were fairly small and may not represent the correct behavior in more challenging problems. 
To ensure the best parameter choice for this new study, we ran the 3-D C5G7 benchmark on a small cluster. 

To identify the best way to use the preconditioner, we looked at different combinations of reasonable values of all of the major parameters the user could choose. 
We used 0.8, 1.0, and 1.2 for the weight in weighted Richardson. 
We used 1, 2, and 3 for both relax count and number of V-cycles. 
Finally, we tried three different subspaces sizes: 50, 100, and 150. 
For all tests, this study used quadruple range (QR) angular quadrature \cite{Jarrell2009} with 8 azimuthal and 6 polar angles per octant (QR 8-6) in the main solver and QR 1-1 inside the preconditioner, $P_1$ scattering expansion, and 56 energy groups. 
The RQI solver used an eigenvalue flux tolerance of $1\times 10^{-3}$ and $k$ tolerance of $1\times 10^{-5}$. 

\begin{figure}
\caption{MGE Parameter Study with 3D-C5G7}
\includegraphics[width=\textwidth]{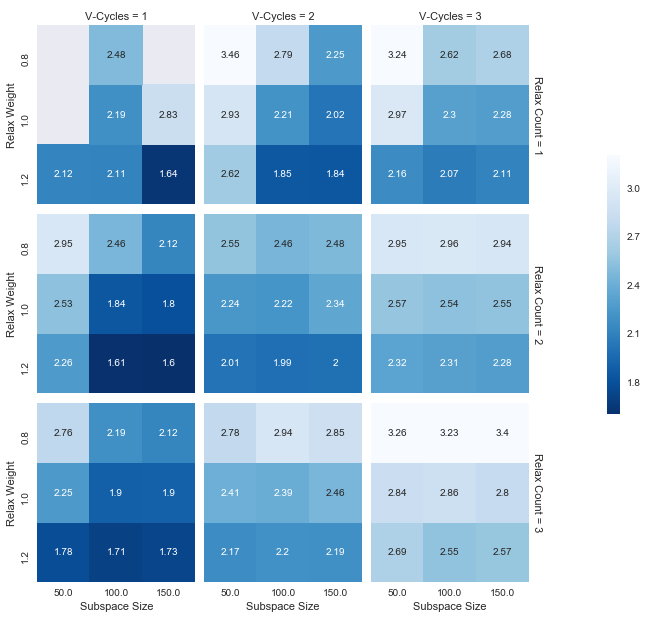}
\label{fig:mge-dataframe}
\centering
\end{figure}

\autoref{fig:mge-dataframe} shows the results of these many combinations of parameters in a compact way. 
Each square contains the runtime (in $1\times 10^{4}$ seconds) for the combination of parameters associated with that square (explained below). 
The lower the number and the darker the shading, the lower the runtime. 
To read the graphic:
\begin{compactitem}
\item Each major column corresponds to the number of V-cycles, where a major column is a set of 3 blocks going across separated from the next set of 3. 
The first column is $v = 1$ and the last is $v = 3$. 
This is indicated across the top.
\item Each major row corresponds to the number of relaxations per level, where a major row is a set of 3 blocks going down separated from the next set of 3. 
The first row is $r = 1$ and the last is $r = 3$. 
This is indicated along the right.
\item Each minor row corresponds to the weight, where a minor rows is each individual block going down. 
The top minor row inside of each major row is $w = 0.8$ and the bottom is $w = 1.2$. 
This is indicated along the left.
\item Each minor column corresponds to the subspace size, where a minor column is each individual block going across. 
The left minor column inside of each major column is subspace = 50 and the right is 150.
This is indicated across the bottom.
\end{compactitem}
For example, to find $w1.2r1v1$ with subspace 150, we go to the first major column (left most, for $v=1$), first major row (top most, for $r=1$), last minor column (right most inside of $v=1$ to get subspace 150), and last minor row (bottom most inside of $r=1$ to get $w=1.2$). This yields a value of $1.64 \times 10^4$ s. 

The general trends from these results indicate that $w1.2r1v1$ or $w1.2r2v1$ with a subspace size of 150 are the best choices. Weight and number of V-cycles had the clearest trends, and subspace size had the least clear trend.

\subsection{PWR 900}
\label{subsec:PWR}
We next studied a large problem that is ``grand challenge'' in scale, a full pressurized water reactor (PWR) 900 core \cite{evans2009d}. This problem used 2 $\times$ 2 spatial cells/pin for 17 $\times$ 17 pins/assembly and 289 assemblies (132 reflector, 159 fuel of varying enrichment). This gave 578 $\times$ 578 $\times$ 700 mesh elements (233,858,800 cells). We used a $P_0$ scattering expansion, an $S_{12}$ angular quadrature, and 44 energy groups--all of which gave 1.73 trillion unknowns. Based on PWR calculations done previously by Evans and Davidson \cite{Evans2010}, $k$ is approximately 1.27. 

\subsubsection{RQI vs.\ PI}
We first did a study on Jaguar \cite{jaguar2010} to better characterize the performance of preconditioned RQI and compare it to preconditioned PI. The preconditioner settings were $w1r3v3$\footnote{Note that scoping calculations showed this problem needed more preconditioning to converge, hence we did not use the just-suggested defaults.}; 4 and 11 sets were used giving 11 and 4 groups per set with 6 and 4 energy grids, respectively. We used 102 $x$-blocks, 100 $y$-blocks, and 10 $z$-blocks. The results are in Table~\ref{table:full PWR}. ``Krylov" and ``Eigenvalue" indicate the number of iterations needed for the respective sovlers. This problem used tolerance = 1e-3, upscattering tolerance = 1e-4, and $k$ tolerance = 1e-3.
\begin{table}[!h]
\caption{PWR 900 Preconditioned Strong Scaling Study with RQI and PI}
\begin{center}
\begin{tabular}{| l | c | c | c | l | c | c | l |}
\hline
Solver & Sets & Cores & $k$ & Krylov & Eigenvalue & Time (m) \\[0.5ex]
\hline
RQI & 4   & 40,800   & 1.269 $\pm$ 1.12e-3 & 76   &  6               & 802.60  \\
PI    & 4   & 40,800   & 1.270 $\pm$ 6.68e-2 & 101 & 10$^{+}$ & 1440.12  \\
RQI & 11 & 112,200 & 1.269 $\pm$ 1.12e-3 & 76   & 6                & 331.43    \\
PI    & 11 & 112,200 & 1.270 $\pm$ 5.09e-2 & 111 & 11$^{+}$ & 480.63     \\
\hline 
\end{tabular}\\
$^{+}$exceeded wall time limit 
\end{center}
\label{table:full PWR}
\end{table}  

Preconditioned RQI was significantly faster than preconditioned PI. In particular, PI was not able to converge before the wall time limit was reached. A true comparison between PI and RQI is difficult because PI never finished the calculations. However, the results show that RQI was much faster and required far fewer Krylov and eigenvalue iterations than PI for this problem. 

We also considered what would happen with a reduced quadrature set, in this case $S_{2}$, inside the preconditioner. With 4 sets, the time with RQI was reduced from 802.60 to 192.48 minutes and the same number of Krylov and RQI iterations were needed. Using the reduced quadrature in MGE had a big time saving impact: reducing the solver time by 76\% without affecting the solution. All further PWR 900 calculations used $S_{2}$ in the preconditioner.

To more conclusively resolve the RQI-to-PI comparison, we did another study, this time on Titan. We used a V-cycle depth of 2, a reduced quadrature in the MGE preconditioner of $S_2$, and tolerances of $1 \times 10^{-3}$. This time we broke the problem up over 112 $\times$ 112 $\times$ 10 partitions (12,544 blocks).

The results using 11 energy sets are given in Table~\ref{tab:PWR all} and illustrate several things. One result is that the MGE preconditioner does not
help PI; when we add it, the calculation slows down enough that it did not
finish within available wall time limits, although convergence does occur
eventually.  The most important result is that preconditioned RQI can be much
faster than PI in general.  Preconditioned RQI was better whether PI was
preconditioned or not. RQI with MGE was more than 10 times faster than PI.
\begin{table}[!h]
\begin{center}
\caption{PWR 900 Comparison of PI and RQI with and without Preconditioning, 11 Energy Sets}
\label{tab:PWR all}
    \begin{tabular}{| c | c | c | c | c | c | c |}
      \hline
      Method & Precond. & N Eigen & N Krylov & $k$ & time (m) \\\hline
      PI   & none     & 149 & 5602 & 1.276 $\pm$ 1.85e-3 & 612.2 \\
      PI   & w1r2v2 & 86   & 946   & 1.275 $\pm$ 1.43e-3 & 720$^+$ \\
      RQI & w1r2v2 & 5     & 70     & 1.268 $\pm$ 1.24e-3 & 54.8 \\
      \hline
    \end{tabular}\\
    $^{+}$exceeded walltime limit\\
  \end{center}
\end{table}

\subsubsection{RQI Strong Scaling}
After demonstrating that RQI can be better than PI for a real problem, we looked at how it performs in a strong scaling study. 
We used RQI with MGE using $w1r2v2$ and 1, 4, 11, and 22 sets. 
The results are given in Table~\ref{tab:PWR rqi strong scaling} where $\text{t}_{\text{perfect}}$ = ($\text{1 set solve time}$ / $\text{\# energy sets}$) and efficiency = ($\text{t}_{\text{perfect}}$ / $\text{t}_{\text{actual}}$). 
A strong scaling study with MGE has been published before, but in that case the V-cycle depth was not fixed.
This meant that increasing energy sets decreased V-cycle depth such that the preconditioner did less work with more sets. 
In this study, the amount of work done by the preconditioner does not vary with the number of energy sets. 
\begin{table}[!h]
\caption{PWR 900 RQI Strong Scaling $w1r2v2$ MGE Preconditioning on Titan}
  \begin{center}
    \begin{tabular}{| c | c | c | c | c | c | c |}
     \hline
      Sets & Domains & N Eigen & N Krylov & time (m) & t$_{\text{perfect}}$ & Efficiency \\\hline
      1   & 12,544    & 5 & 70 & 407.8 & 407.8 & 1.000 \\
      4   & 50, 176   & 5 & 70 & 123.4 & 102.0 & 0.826\\
      11  & 137,984 & 5 & 70 & 54.8  & 37.1  & 0.676\\
      22  & 275,968 & 5 & 70 & 39.6  & 18.5  & 0.468\\
      \hline
    \end{tabular}
    \label{tab:PWR rqi strong scaling}
  \end{center}
\end{table}
\begin{figure}
\caption{PWR 900 RQI Strong Scaling $w1r2v2$ MGE Preconditioning on Titan}
\includegraphics[scale=.75]{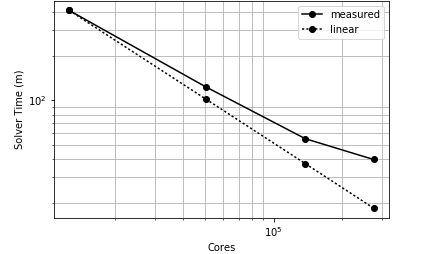}
\centering
\label{fig:pwr-scaling}
\end{figure}

The scaling, which is plotted in Figure~\ref{fig:pwr-scaling}, compares quite well to previous scaling studies for Denovo. A fixed source (i.e.\ MG Krylov only) problem with a similar mesh and 44 groups scaled from 4,320 domains to 190,080 domains with an efficiency of 0.64 \cite{Slaybaugh2011}. That this problem performed similarly shows that adding RQI and the MGE preconditioner as solvers does not degrade the strong scaling achieved using the MG Krylov solver only. It is promising that the new solver system does not degrade scaling and that a problem for which RQI is decisively faster than PI is one for which this work was designed.

\subsection{BW-1484}
Finally, we investigated performance of the new solvers with a model of the Babcock and Wilcox 1484 reactor (BW1484) \cite{bw1484}. We used this problem to more deeply investigate which solvers work well for real problems and under what conditions. We looked at scaling in energy for preconditioned and unpreconditioned RQI with Arnoldi using both upscatter and regular partitioning. 


We started by performing a medium-sized scaling study on a small cluster. This problem used 2 $\times$ 2 spatial cells/pin with 60 $\times$ 60 pins for a total of 120 $\times$ 120 $\times$ 44 mesh elements (633,600 cells) split up across 8 $x$-blocks,16 $y$-blocks, and 1 $z$-block. 
We used a $P_0$ scattering expansion; a QR 8-6 angular quadrature; and 56 energy groups--which gave approximately 425 million unknowns. 
We used a $k$ tolerance and flux tolerance of $1 \times 10^{-5}$.
In the preconditioned cases, a QR 1-1 quadrature was used in the multilevel preconditioner with $w1.2r2v1$. 
The study was performed with 128, 256, 512, and 1024 cores and 1, 2, 4, and 8 energy sets, respectively. 

\begin{figure}
\caption{BW1484 Energy Small Energy Scaling Study}
\includegraphics[scale=1]{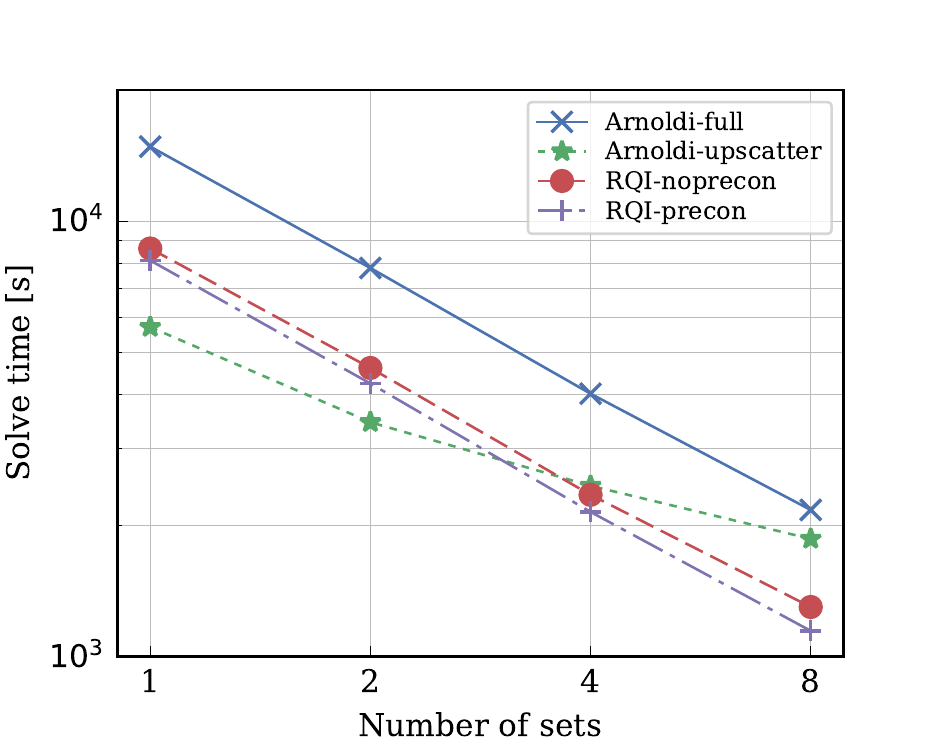}
\label{fig:small}
\centering
\end{figure}
Figure \ref{fig:small} shows the results from this scaling study. We can see favorable scaling with both preconditioned and unpreconditioned RQI. Preconditioned RQI was faster than unpreconditioned. Both RQI results were faster than Arnoldi with full partitioning (which scaled just as well, but was the slowest). Arnoldi with upscatter partitioning performed the best at low core count, but shows poor scaling performance and, after two energy sets, ended in third place with respect to time. 

Finally, we used Titan to look at scaling performance with a larger version of this problem with more cores. Here the 1-set case was run with 3,136 cores (56 $x$-blocks and 56 $y$-blocks). A parallel-efficiency modeling tool developed at Oak Ridge was used to choose the optimal number of $z$-blocks for each method, which in the case of RQI with MGE preconditioning was determined to be 22. For the other three methods, using 11 $z$-blocks was determined to be optimal. We then scaled up to 2, 4, and 8 sets using 6,272, 12,544, and 25,088 cores, respectively. The upscatter and $k$ tolerances were $1 \times 10^{-5}$. We used 8 mesh cells/pin for a total of 336 $\times$ 336 $\times$ 154 mesh elements (17,385,984 cells). We again used $P_0$ and QR 8-6, with $w1.2r2v1$ and QR 1-1 inside the preconditioner.

\begin{figure}
\caption{BW1484  Large Energy Scaling Study on Titan}
\includegraphics[scale=1]{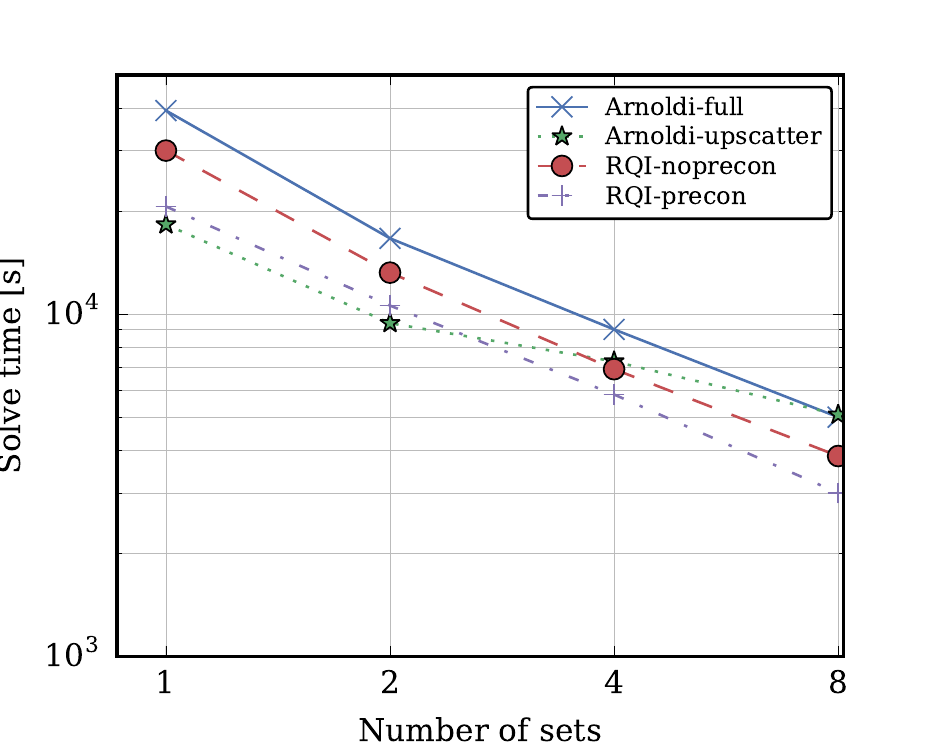}
\label{fig:titan}
\centering
\end{figure}
The results can be seen in Figure~\ref{fig:titan}. The same trends were seen on Titan as in the small study. Arnoldi with full partitioning scaled well, but was the slowest. RQI with and without preconditioning scaled similarly, with the preconditioned version being faster. Arnoldi with upscatter partitioning started out as the fastest, but scales the least well. Thus, at 4 and 8 energy sets it was much slower. 

We can see that RQI gives the same scaling in energy relative to Arnoldi with full partitioning. We expected this behavior as they handle energy sets the same way. Therefore, RQI is the clear winner in terms of runtime when energy groups are handled the same way. 

Furthermore, the Arnoldi behavior makes sense given the way the block Krylov solver works and the way partitioning works. With full partitioning, all of the groups are divided into energy sets and solved with a Krylov solver. Thus, adding more sets increases the parallelization in energy. With only the upscattering partitioned, the downscatter only groups, which are half of the groups in this case, are always solved by every energy set solved with GS. Thus, using more sets does not result in additional parallelization in energy after the upscatter groups are fully decomposed and we see the scaling stagnate. 

Another way to look at the full vs.\ upscatter partitioning helps clarify why upscatter partitioning is faster with one energy set. With upscatter partitioning, each downscatter-only group is solved with GS--which is a one-group-sized Krylov solve and converges the group on the first multigroup iteration. That's tough to beat. With full partitioning, we're putting the downscattering groups with the upscattering groups into an all-groups-sized Krylov solve. We're now using an iteration vector that is much larger and needs all groups to converge for any group to converge. We expect that in many cases this will slow down converging the downscatter groups. 

These results lead us to conclude that for real, challenging problems, Arnoldi with upscatter partitioning should be used when energy parallelization is not needed. For problems that are large enough to need energy parallelization, RQI with the MGE preconditioner should be selected. This makes sense given the structure of the solvers.

\section{Conclusions}
\label{sec:conclusions}
The goal of this research was to accelerate transport calculations for hard problems with methods that can take full advantage of modern leadership-class computers, facilitating the design of better nuclear systems. Three complimentary methods were implemented that accomplish this goal. 

At the outset of this work, Denovo could be decomposed in five of the six dimensions of phase space over which the steady-state transport equation is solved in a way that restricted it to about 20,000 cores for a problem with 500 million cells. The original suite of solvers included GS as the multigroup solver and PI with inner GS iteration as the eigenvalue solver. These solvers have some significant limitations in many cases of interest. 

The new MG Krylov solver converges more quickly than GS and enables energy decomposition such that Denovo can scale to hundreds of thousands of cores. The new MGE preconditioner reduces iteration count for many problem types and takes advantage of the new energy decomposition such that it can scale very efficiently. These two tools are useful on their own, but together they allow the RQI eigenvalue solver to work.

The real motivation of this work was to add RQI, which should converge in fewer iterations and less time than PI, and possibly Arnoldi, for large and challenging problems. RQI creates shifted systems that would not be tractable without the MG Krylov solver. It also creates ill-conditioned matrices that cannot converge without the MGE preconditioner. Using these methods, RQI converged in fewer iterations and in less time than both PI and Arnoldi for large problems on large core counts. 

The methods added in this research accelerated Denovo in multiple ways. This acceleration helps enable the solution of today's ``grand challenge'' problems. It is hoped that improved methods will lead to improved reactor designs and systems, and that the frontier of computational challenges will be moved forward.

%

This set of solvers has never been previously combined and demonstrated, nor has RQI been applied to the transport equation before this line of research. This is likely because it takes $O(n^{3})$ operations for full, dense matrices, and without parallelization in energy it could be prohibitively expensive \cite{Stewart2001}. In the past, adding a shift to make the scattering matrix energy-block dense was difficult to handle. The system would have been solved with GS, which would have been restrictively slow. The MG Krylov algorithm has enabled energy parallelization and made the calculation of the eigenvector tractable. 

The combination of the MG block Krylov solver, RQI, and the MGE preconditioner performed very well for large problems. It also scaled well in energy, enabling calculations on large core counts. We expect that this collection of solvers can be highly valuable for very challenging calculations needed for detailed multiphysics, perturbation studies, and physics investigations. 



\pagebreak
\section*{Acknowledgments}

This research used resources of the Oak Ridge Leadership Computing Facility at the Oak Ridge National Laboratory, which is supported by the Office of Science of the U.S. Department of Energy under Contract No. DE-AC0500OR22725. Additional thanks to the Rickover Fellowship Program in Nuclear Engineering sponsored by Naval Reactors Division of the U.S.\ Department of Energy. This fellowship sponsored the work from which this work is derived.

\pagebreak

\bibliographystyle{nse}
\bibliography{RQI_MGE}

\end{document}